\documentclass[12pt]{article}

\usepackage{amsmath, epsfig, cite}
\usepackage{amssymb}
\usepackage{amsfonts}
\usepackage{latexsym}
\usepackage{float}
\usepackage{color}

\newtheorem{thm}{Theorem}[section]

\newtheorem{lem}[thm]{Lemma}

\numberwithin{equation}{section}

\newcommand{\qed}{{\hfill$\square$}\medskip}
\UseRawInputEncoding
\begin{document}

\begin{center}
{\Large\bf Supercongruences involving  Motzkin numbers\\[7pt]
 and central trinomial coefficients}
\end{center}

\vskip 2mm \centerline{Ji-Cai Liu}
\begin{center}
{\footnotesize Department of Mathematics, Wenzhou University, Wenzhou 325035, PR China\\
{\tt jcliu2016@gmail.com } \\[10pt]
}
\end{center}


\vskip 0.7cm \noindent{\bf Abstract.}
Let $M_n$ and $T_n$ denote the $n$th Motzkin number and the $n$th central trinomial coefficient respectively. We prove that for any prime $p\ge 5$,
\begin{align*}
&\sum_{k=0}^{p-1}M_k^2\equiv \left(\frac{p}{3}\right)\left(2-6p\right)\pmod{p^2},\\
&\sum_{k=0}^{p-1}kM_k^2\equiv \left(\frac{p}{3}\right)\left(9p-1\right)\pmod{p^2},\\
&\sum_{k=0}^{p-1}T_kM_k\equiv \frac{4}{3}\left(\frac{p}{3}\right)+\frac{p}{6}\left(1-9\left(\frac{p}{3}\right)\right)\pmod{p^2},
\end{align*}
where $\left(-\right)$ is the Legendre symbol. These results confirm three 12-year-old supercongruence conjectures of Z.-W. Sun.

\vskip 3mm \noindent {\it Keywords}: Motzkin numbers; central trinomial coefficients; supercongruences

\vskip 2mm
\noindent{\it MR Subject Classifications}: 11B65, 11A07, 05A19
\section{Introduction}
The Motzkin numbers $\{M_n\}_{n=0}^{\infty}=1, 1, 2, 4, 9, 21, 51, 127,\cdots$ first appeared in \cite{motzkin-bams-1948} in a circle chording setting, which count the number of ways of connecting a subset of $n$ points on a circle by nonintersecting chords.
The Motzkin number $M_n$ also counts the number of lattice paths on the upper right quadrant of a grid from $(0,0)$ to $(n,0)$ if one is allowed to move by using only steps $(1,1),(1,0)$ and $(1,-1)$ but forbidden from dipping below the $y=0$ axis.

The Motzkin numbers are named after Theodore Motzkin and naturally appear in various combinatorial objects. Fourteen different manifestations of Motzkin numbers in different branches of mathematics were enumerated by Donaghey and Shapiro \cite{ds-jcta-1977} in their survey of Motzkin numbers. The interested readers may refer to \cite{stanley-b-1999} for further information on Motzkin numbers.

The famous Catalan numbers $C_n={2n\choose n}/(n+1)$ are closely related to Motzkin numbers.
The Motzkin numbers can be expressed in terms of Catalan numbers:
\begin{align*}
M_n=\sum_{k=0}^{\lfloor n/2\rfloor}{n\choose 2k}C_k,
\end{align*}
and inversely,
\begin{align*}
C_{n+1}=\sum_{k=0}^n{n\choose k}M_k,
\end{align*}
where $\lfloor x \rfloor$ denotes the integral part of real $x$.

Another sequence closely related to Motzkin numbers is the central trinomial coefficients. For $n\in \mathbb{N}$, the central trinomial coefficient $T_n$ is given by the constant term in the expansion of $(1+x+x^{-1})^2$, which can be expressed in terms of binomial coefficients:
\begin{align*}
T_n=\sum_{k=0}^{\lfloor n/2\rfloor}{n\choose 2k}{2k\choose k}.
\end{align*}
We remark that the central trinomial coefficient $T_n$ counts the number of lattice paths from $(0,0)$ to $(n,0)$ if one is allowed to move by using only steps $(1,1),(1,0)$ and $(1,-1)$.

Although Catalan numbers, Motzkin numbers and central trinomial coefficients naturally arise in combinatorics, they also possess rich arithmetic properties.

Throughout the paper, let $p\ge 5$ be a prime. Sun and Tauraso \cite{st-ijnt-2011} showed that
\begin{align*}
\sum_{k=0}^{p-1}C_k\equiv\frac{3}{2}\left(\frac{p}{3}\right)- \frac{1}{2}\pmod{p^2},
\end{align*}
where $\left(-\right)$ denotes the Legendre symbol.

In 2014, Sun \cite{sunzw-scm-2014} proved that
\begin{align*}
&\sum_{k=0}^{p-1}T_k^2\equiv \left(\frac{-1}{p}\right)\pmod{p},\\
&\sum_{k=0}^{p-1}T_kM_k\equiv \frac{4}{3}\left(\frac{p}{3}\right)\pmod{p}.
\end{align*}
By the telescoping method developed by Chen, Hou and Mu \cite{chm-jcam-2006}, Sun \cite{sunzw-aam-2022} recently established the following interesting supercongruence:
\begin{align}
\sum_{k=0}^{p-1}(2k+1)M_k^2\equiv 12p\left(\frac{p}{3}\right)\pmod{p^2}.\label{a-1}
\end{align}

For more arithmetic properties of Catalan numbers, Motzkin numbers and central trinomial coefficients, we refer the interested readers to \cite{cp-hjm-2014,cw-pams-2022,ds-jnt-2006,ely-ejc-2008,
km-ejc-2018,liu-cmr-2020,liu-rj-2022,ps-dm-2006,tauraso-aam-2012}.

The motivation of the paper is the following three conjectural supercongruences of Sun \cite[Conjecture 1.1]{sunzw-scm-2014}:
\begin{align}
&\sum_{k=0}^{p-1}M_k^2\equiv \left(\frac{p}{3}\right)\left(2-6p\right)\pmod{p^2},\label{a-2}\\
&\sum_{k=0}^{p-1}kM_k^2\equiv \left(\frac{p}{3}\right)\left(9p-1\right)\pmod{p^2},\label{a-3}\\
&\sum_{k=0}^{p-1}T_kM_k\equiv \frac{4}{3}\left(\frac{p}{3}\right)+\frac{p}{6}\left(1-9\left(\frac{p}{3}\right)\right)\pmod{p^2}.
\label{a-4}
\end{align}
``The three supercongruences look curious and challenging", as described by Sun \cite{sunzw-aam-2022} in his recent paper. Although the three supercongruence conjectures were officially announced by Sun \cite{sunzw-scm-2014} in 2014, they first appeared in arXiv version of Sun's paper in 2010 (see  https://arxiv.org/abs/1008.3887), and have a history of 12 years.

In this paper, we aim to prove supercongruences \eqref{a-2}--\eqref{a-4}.
\begin{thm}\label{t-1}
The supercongruences \eqref{a-2}--\eqref{a-4} are true.
\end{thm}

We remark that both Motzkin numbers and central trinomial coefficients have many different
expressions (see A001006 and A002426 in \cite{sloane-1964}). The following two expressions will be used in the proof of Theorem \ref{t-1}:
\begin{align}
&M_n=\sum_{k=0}^n (-1)^{n+k}{n\choose k}C_{k+1},\label{a-5}\\
&T_n=\sum_{k=0}^n (-1)^{n+k}{n\choose k}{2k\choose k}.\label{a-6}
\end{align}

The rest of the paper is organized as follows. Section 2 is devoted to proving some preliminary results. In Section 3, we establish three congruences for triple sums which paly an important role in the proof of Theorem \ref{t-1}. We prove \eqref{a-2} and \eqref{a-3} in Section 4, and prove \eqref{a-4} in Section 5.

\section{Preliminaries}
In the proof of Theorem \ref{t-1}, we will frequently use Wolstenholme's theorem \cite{wols-1862}:
\begin{align}
{2p-1\choose p-1}\equiv 1\pmod{p^3},\label{b-1}
\end{align}
which is equivalent to
\begin{align}
{2p\choose p}\equiv 2\pmod{p^3},\label{b-2}
\end{align}
and Lehmer's congruences \cite{lehmer-am-1938}:
\begin{align}
&H_{\lfloor p/6\rfloor}\equiv H_{\lfloor 5p/6\rfloor}\equiv -2q_p(2)-\frac{3}{2}q_p(3)\pmod{p},\label{bL-1}\\
&H_{\lfloor p/3\rfloor}\equiv H_{\lfloor 2p/3\rfloor}\equiv -\frac{3}{2}q_p(3)\pmod{p},\label{bL-2}\\
&H_{\lfloor p/2\rfloor}\equiv -2q_p(2)\pmod{p},\label{bL-3}
\end{align}
where $q_p(a)$ denotes the Fermat quotient $(a^{p-1}-1)/p$.

In addition, we require some congruences related to central binomial coefficients and Catalan numbers.
\begin{lem}
For any prime $p\ge 5$, we have
\begin{align}
&\sum_{k=0}^{p-1}{2k\choose k}\equiv \left(\frac{p}{3}\right)\pmod{p^2},\label{b-3}\\
&\sum_{k=0}^{p-1}C_k\equiv\frac{3}{2}\left(\frac{p}{3}\right)- \frac{1}{2}\pmod{p^2},\label{b-4}\\
&\sum_{k=1}^{p-1}\frac{{2k\choose k}}{k}\equiv 0\pmod{p^2},\label{bnew-4}\\
&\sum_{k=0}^{p-1}k{2k\choose k}\equiv\frac{2}{3}\left(p-\left(\frac{p}{3}\right)\right)\pmod{p^2},\label{b-5}\\
&\sum_{k=0}^{p-2}\frac{{2k\choose k}}{(k+1)^2}\equiv 3\left(\frac{p}{3}\right)+1\pmod{p}.
\label{b-18}
\end{align}
\end{lem}
{\noindent\it Proof.}
Congruences \eqref{b-3}--\eqref{bnew-4} were proved by Sun and Tauraso (see \cite[(1.7) and (1.9)]{st-ijnt-2011} and \cite[Theorem 1.3]{st-aam-2010}). It is easily proved by induction on $n$ that
\begin{align}
\sum_{k=0}^{n-1}(3k+2){2k\choose k}=n{2n\choose n}.\label{b-6}
\end{align}
Letting $n=p$ in \eqref{b-6} and using \eqref{b-2} gives
\begin{align}
\sum_{k=0}^{p-1}(3k+2){2k\choose k}\equiv 2p\pmod{p^4}.\label{b-7}
\end{align}
Then the proof of \eqref{b-5} follows from \eqref{b-3} and \eqref{b-7}.

From the following identity:
\begin{align*}
{2k\choose k+1}=\frac{1}{2}{2k+2\choose k+1}-{2k\choose k},
\end{align*}
we deduce that
\begin{align}
\sum_{k=0}^{p-2}\frac{{2k\choose k+1}}{k+1}&=\frac{1}{2}\sum_{k=1}^{p-1}\frac{{2k\choose k}}{k}
-\sum_{k=0}^{p-2}C_k\notag\\[7pt]
&\equiv \frac{1}{2}-\frac{3}{2}\left(\frac{p}{3}\right)+C_{p-1}\notag\\[7pt]
&\equiv -\frac{1}{2}-\frac{3}{2}\left(\frac{p}{3}\right)\pmod{p},\label{bnew-9}
\end{align}
where we have used \eqref{b-4} and \eqref{bnew-4}.

By \eqref{b-4}, \eqref{bnew-9} and the identity $C_k={2k\choose k}-{2k\choose k+1}$, we have
\begin{align*}
\sum_{k=0}^{p-2}\frac{{2k\choose k}}{(k+1)^2}&=\sum_{k=0}^{p-2}C_k-\sum_{k=0}^{p-2}\frac{{2k\choose k+1}}{k+1}\\[7pt]
&\equiv 3\left(\frac{p}{3}\right)+1\pmod{p},
\end{align*}
as desired.
\qed

We also need some congruences involving harmonic numbers $H_n=\sum_{k=1}^n 1/k$.
\begin{lem}
For any prime $p\ge 5$, we have
\begin{align}
&\sum_{k=0}^{p-1}{2k\choose k}H_{k}\equiv -\left(\frac{p}{3}\right)q_p(3) \pmod{p},\label{b-8}\\
&\sum_{k=0}^{p-1}k{2k\choose k}H_{k}\equiv \frac{1}{3}\left(\frac{p}{3}\right)\left(\left(\frac{p}{3}\right)+2q_p(3)-1\right)\pmod{p},\label{bnew-3}\\
&\sum_{k=0}^{p-1}C_{k}H_k\equiv -\frac{3}{2}\left(\frac{p}{3}\right)q_p(3)\pmod{p}.\label{b-9}
\end{align}
\end{lem}
{\noindent\it Proof.}
We remark that \eqref{b-8} was already mentioned in \cite[page 528]{ms-ijnt-2016}.
By using the software package {\tt Sigma} developed by Schneider \cite{schneider-slc-2007}, we
can automatically discover and prove the following three combinatorial identities:
\begin{align}
&\sum_{k=0}^{n}(-4)^k {n\choose k}H_{k}=(-3)^n\left(H_n-\sum_{k=1}^n\frac{1}{k(-3)^k}\right),\label{b-10}\\[7pt]
&\sum_{k=0}^{n}(-4)^k k{n\choose k}H_{k}=\frac{1-(-3)^n}{3}+\frac{4n(-3)^n}{3}\left(H_n-\sum_{k=1}^n
\frac{1}{k(-3)^k}\right),\label{bnew-1}
\end{align}
and
\begin{align}
&\sum_{k=0}^n\frac{(-4)^k}{k+1}{n\choose k}H_k\notag\\
&=\frac{1}{4(n+1)}\left((-1+3(-3)^n)H_n-3(-3)^n
\sum_{k=1}^n\frac{1}{k(-3)^k}+\sum_{k=1}^n\frac{(-3)^k}{k}\right).\label{b-11}
\end{align}
By \eqref{b-10}, we can rewrite \eqref{bnew-1} and \eqref{b-11} as
\begin{align}
\sum_{k=0}^{n}(-4)^k k{n\choose k}H_{k}=\frac{1-(-3)^n}{3}+\frac{4n}{3}\sum_{k=0}^{n}(-4)^k {n\choose k}H_{k},\label{bnew-2}
\end{align}
and
\begin{align}
&\sum_{k=0}^n\frac{(-4)^k}{k+1}{n\choose k}H_k\notag\\[7pt]
&=\frac{1}{4(n+1)}\left(\frac{3((-3)^n+1)}{(-3)^n}\sum_{k=0}^n(-4)^k{n\choose k}H_k\right)\notag\\[7pt]
&+\frac{1}{4(n+1)}\left(3\sum_{k=1}^n\frac{1}{k(-3)^k}+\sum_{k=1}^n\frac{(-3)^k}{k}-4H_n\right).\label{b-12}
\end{align}
Using Fermat's little theorem, we obtain
\begin{align*}
&\sum_{k=1}^{(p-1)/2}\frac{(-3)^k}{k}+3\sum_{k=1}^{(p-1)/2}\frac{1}{k(-3)^{k}}\\[7pt]
&\equiv \sum_{k=1}^{(p-1)/2}\frac{(-3)^k}{k}+\sum_{k=1}^{(p-1)/2}\frac{(-3)^{p-k}}{p-k}\\[7pt]
&=\sum_{k=1}^{p-1}\frac{(-3)^k}{k}\pmod{p}.
\end{align*}
From Granville's congruence \cite{granville-integer-2004}:
\begin{align*}
\sum_{j=1}^{p-1}\frac{x^j}{j}\equiv \frac{1-x^p+(x-1)^p}{p}\pmod{p},
\end{align*}
we deduce that
\begin{align*}
\sum_{k=1}^{p-1}\frac{(-3)^k}{k}\equiv 3q_p(3)-8q_p(2)\pmod{p},
\end{align*}
and so
\begin{align}
\sum_{k=1}^{(p-1)/2}\frac{(-3)^k}{k}+3\sum_{k=1}^{(p-1)/2}\frac{1}{k(-3)^{k}}
\equiv 3q_p(3)-8q_p(2)\pmod{p}.\label{b-13}
\end{align}
Finally, letting $n=(p-1)/2$ in \eqref{bnew-2}--\eqref{b-12} and using \eqref{bL-3}, \eqref{b-8}, \eqref{b-13}
and the facts that
\begin{align*}
&(-3)^{(p-1)/2}\equiv \left(\frac{p}{3}\right)\pmod{p},\\[7pt]
&{(p-1)/2\choose k}\equiv \frac{{2k\choose k}}{(-4)^k}\pmod{p},
\end{align*}
we arrive at the desired congruences \eqref{bnew-3} and \eqref{b-9}.
\qed

For an assertion $A$, we set
\begin{align*}
[A]=\begin{cases}
1\quad&\text{if $A$ holds,}\\
0\quad&\text{otherwise.}
\end{cases}
\end{align*}
The following two known congruences paly an important role in the proof of Theorem \ref{t-1} (see \cite[Theorem 1.2]{ps-dm-2006}).
\begin{lem}\label{lem-3}
Let $p\ge 5$ be a prime. For $1\le k \le p$, we have
\begin{align}
&\sum_{i=0}^{p-1}{2i\choose i+k}\equiv \left(\frac{p-k}{3}\right)\pmod{p},\label{b-14}\\[7pt]
&\sum_{i=1}^{p-1}\frac{{2i\choose i+k}}{i}\equiv\frac{\alpha(k)-1}{k}\pmod{p},\label{b-15}
\end{align}
where $\alpha(k)=2(-1)^k+3[3\mid p-k])$.
\end{lem}

Based on Lemma \ref{lem-3}, we establish the following result which seems to be crude but useful in
the proof of Theorem \ref{t-1}.
\begin{lem}\label{lem-4}
Let $p\ge 5$ be a prime. For $1\le k \le p$, we have
\begin{align}
(-1)^k\sum_{i=1}^{p-2}\frac{{2i\choose i+k}}{i+1}
\equiv 2k+\frac{3}{2}\left(\frac{p}{3}\right)-\frac{3}{2}+\sum_{i=1}^{k-1}\beta(i)
-k\sum_{i=1}^{k-1}\frac{\beta(i)+2}{i}\pmod{p},\label{b-16}
\end{align}
where $\beta(i)=(-1)^i(3[3\mid p-i]-1)$.
\end{lem}

{\noindent\it Proof.}
By Pascal's formula ${n\choose m}={n-1\choose m}+{n-1\choose m-1}$, we have
\begin{align*}
{2i+2\choose i+1+k}&={2i+1\choose i+1+k}+{2i+1\choose i+k}\\[7pt]
&={2i \choose i+1+k}+2{2i\choose i+k}+{2i\choose i+k-1}.
\end{align*}
It follows that
\begin{align}
\sum_{i=0}^{p-2}\frac{{2i \choose i+1+k}}{i+1}+2\sum_{i=0}^{p-2}\frac{{2i\choose i+k}}{i+1}+\sum_{i=0}^{p-2}\frac{{2i\choose i+k-1}}{i+1}
=\sum_{i=1}^{p-1}\frac{{2i\choose i+k}}{i}.\label{bnew-5}
\end{align}

Let
\begin{align*}
&f(k)=(-1)^k\sum_{i=0}^{p-2}\frac{{2i\choose i+k}}{i+1},\\[7pt]
&g(k)=(-1)^{k+1}\sum_{i=1}^{p-1}\frac{{2i\choose i+k}}{i},\\[7pt]
&F(k)=f(k+1)-f(k).
\end{align*}
We rewrite \eqref{bnew-5} as
\begin{align}
F(k)-F(k-1)=g(k).\label{bnew-6}
\end{align}
From \eqref{bnew-6}, we deduce that
\begin{align*}
F(k)=f(k+1)-f(k)=\sum_{j=1}^k g(j)+F(0),
\end{align*}
and so
\begin{align}
f(k)&=\sum_{l=1}^{k-1}\sum_{j=1}^l g(j)+kF(0)+f(0)\notag\\[7pt]
&=\sum_{j=1}^{k-1}(k-j) g(j)+kf(1)+(1-k)f(0).\label{bnew-7}
\end{align}

By ${2i\choose i+1}={2i\choose i}-{2i\choose i}/(i+1)$, we have
\begin{align}
&kf(1)+(1-k)f(0)\notag\\[7pt]
&=-k\sum_{i=0}^{p-2}\frac{{2i\choose i+1}}{i+1}+(1-k)\sum_{i=0}^{p-2}\frac{{2i\choose i}}{i+1}\notag\\[7pt]
&=(1-2k)\sum_{i=0}^{p-2}C_i+k\sum_{i=0}^{p-2}\frac{{2i\choose i}}{(i+1)^2}\notag\\[7pt]
&\equiv \frac{3}{2}\left(\frac{p}{3}\right)+\frac{1}{2}\pmod{p},\label{bnew-8}
\end{align}
where we have used \eqref{b-4} and \eqref{b-18} in the last step.

Combining \eqref{b-15}, \eqref{bnew-7} and \eqref{bnew-8} gives
\begin{align*}
&(-1)^k\sum_{i=1}^{p-2}\frac{{2i\choose i+k}}{i+1}\\[7pt]
&\equiv -k\sum_{i=1}^{k-1}
\frac{(-1)^i(3[3\mid p-i]-1)+2}{i}+\sum_{i=1}^{k-1}(-1)^i(3[3\mid p-i]-1)\\[7pt]
&+2k+\frac{3}{2}\left(\frac{p}{3}\right)-\frac{3}{2}\\[7pt]
&=2k+\frac{3}{2}\left(\frac{p}{3}\right)-\frac{3}{2}+\sum_{i=1}^{k-1}\beta(i)
-k\sum_{i=1}^{k-1}\frac{\beta(i)+2}{i}\pmod{p},
\end{align*}
as desired.
\qed

\section{Three key triple sums}
The main idea in the proof of Theorem \ref{t-1} is to translate the left-hand sides of
\eqref{a-2} and \eqref{a-4} into three tripe sums, which can be determined modulo $p$ by using
Lemmas \ref{lem-3} and \ref{lem-4}. The three congruences for triple sums are stated as follows.
\begin{lem}
For any prime $p\ge 5$, we have
\begin{align}
\sum_{k=1}^{p-2}\sum_{j=1}^{p-2}\sum_{i=1}^{p-1}\frac{(-1)^k}{i(j+1)}{2j\choose j+k}{2i\choose i+k}\equiv \left(\frac{p}{3}\right)\left(1-q_p(3)\right)-1\pmod{p}.\label{c-1}
\end{align}
\end{lem}
{\noindent\it Proof.}
Substituting \eqref{b-15} and \eqref{b-16} into the left-hand side of \eqref{c-1} gives
\begin{align}
&\sum_{k=1}^{p-2}(-1)^k\sum_{j=1}^{p-2}\frac{1}{j+1}{2j\choose j+k}\sum_{i=1}^{p-1}\frac{1}{i}{2i\choose i+k}\notag\\[7pt]
&\equiv \sum_{k=1}^{p-2} \left(1-\alpha(k)\right) \sum_{i=1}^{k-1}\frac{\beta(i)+2}{i}
+\sum_{k=1}^{p-2} \frac{\alpha(k)-1}{k}\sum_{i=1}^{k-1}\beta(i)\notag\\[7pt]
&+2\sum_{k=1}^{p-2} \left(\alpha(k)-1\right)
+\left(\frac{3}{2}\left(\frac{p}{3}\right)-\frac{3}{2}\right)\sum_{k=1}^{p-2} \frac{\alpha(k)-1}{k}\pmod{p}.\label{c-2}
\end{align}

We shall only prove the case $p\equiv 1\pmod{3}$ of \eqref{c-1}. The proof of the case $p\equiv 2\pmod{3}$ runs analogously, and we omit the details.

Suppose that $p\equiv 1\pmod{3}$. Let $n=\lfloor p/6\rfloor$. Then $n\equiv -1/6\pmod{p}$.
Since the sequences $\{\alpha(k)\}_{k\in \mathbb{N}}$ and $\{\beta(k)\}_{k\in \mathbb{N}}$ both have a period of $6$, it is easy to check that
\begin{align*}
\Omega(k)=\sum_{i=1}^{k} \left(\alpha(i)-1\right)=
\begin{cases}
0\quad&\text{for $k\equiv 0,1\pmod{6}$},\\
1\quad&\text{for $k\equiv 2\pmod{6}$},\\
-2\quad&\text{for $k\equiv 3\pmod{6}$},\\
2\quad&\text{for $k\equiv 4\pmod{6}$},\\
-1\quad&\text{for $k\equiv 5\pmod{6}$},
\end{cases}
\end{align*}
and
\begin{align*}
\left(\alpha(k)-1\right)\sum_{i=1}^{k-1}\beta(i)=
\begin{cases}
1\quad&\text{for $k\equiv 0\pmod{6}$},\\
0\quad&\text{for $k\equiv 1,5\pmod{6}$},\\
-2\quad&\text{for $k\equiv 2\pmod{6}$},\\
9\quad&\text{for $k\equiv 3\pmod{6}$},\\
-8\quad&\text{for $k\equiv 4\pmod{6}$}.
\end{cases}
\end{align*}
It follows that
\begin{align}
\Omega(p-2)=-1,\label{c-3}
\end{align}
and
\begin{align}
\sum_{k=1}^{p-2} \frac{\alpha(k)-1}{k}\sum_{i=1}^{k-1}\beta(i)
=\sum_{k=1}^{n-1}\frac{1}{6k}-2\sum_{k=1}^{n}\frac{1}{6k-4}
+9\sum_{k=1}^{n}\frac{1}{6k-3}-8\sum_{k=1}^{n}\frac{1}{6k-2}.\label{c-4}
\end{align}

Furthermore, we have
\begin{align*}
&\sum_{k=1}^{p-2} \left(1-\alpha(k)\right) \sum_{i=1}^{k-1}\frac{\beta(i)+2}{i}\\[7pt]
&=-\sum_{i=1}^{p-3}\frac{\beta(i)+2}{i} \sum_{k=i+1}^{p-2} \left(\alpha(k)-1\right)\\[7pt]
&=-\sum_{i=1}^{p-3}\frac{\left(\beta(i)+2\right)\left(\Omega(p-2)-\Omega(i)\right)}{i}\\[7pt]
&=\sum_{i=1}^{p-3}\frac{\left(\beta(i)+2\right)\left(1+\Omega(i)\right)}{i},
\end{align*}
where we have used \eqref{c-3} in the last step. It is also easy to check that
\begin{align*}
\left(\beta(i)+2\right)\left(1+\Omega(i)\right)=
\begin{cases}
1\quad&\text{for $i\equiv 0\pmod{6}$},\\
0\quad&\text{for $i\equiv 1,5\pmod{6}$},\\
2\quad&\text{for $i\equiv 2\pmod{6}$},\\
-3\quad&\text{for $i\equiv 3\pmod{6}$},\\
12\quad&\text{for $i\equiv 4\pmod{6}$},
\end{cases}
\end{align*}
an so
\begin{align}
&\sum_{k=1}^{p-2} \left(1-\alpha(k)\right) \sum_{i=1}^{k-1}\frac{\beta(i)+2}{i}\notag\\[7pt]
&=\sum_{k=1}^{n-1}\frac{1}{6k}+12\sum_{k=1}^{n}\frac{1}{6k-2}-3\sum_{k=1}^{n}\frac{1}{6k-3}
+2\sum_{k=1}^{n}\frac{1}{6k-4}.\label{c-5}
\end{align}
Substituting \eqref{c-3}--\eqref{c-5} into the right-hand side of \eqref{c-2} gives
\begin{align}
&\sum_{k=1}^{p-2}\sum_{j=1}^{p-2}\sum_{i=1}^{p-1}\frac{(-1)^k}{i(j+1)}{2j\choose j+k}{2i\choose i+k}\notag\\[7pt]
&\equiv 2\sum_{k=1}^{n-1}\frac{1}{6k}+4\sum_{k=1}^{n}\frac{1}{6k-2}+6\sum_{k=1}^{n}\frac{1}{6k-3}-2\notag\\[7pt]
&\equiv \frac{1}{3}\sum_{k=1}^n\frac{1}{k}+\frac{2}{3}\sum_{k=1}^n\frac{1}{k+2n}+\sum_{k=1}^n
\frac{1}{k+3n}\notag\\[7pt]
&=\frac{1}{3}H_n-\frac{2}{3}H_{2n}-\frac{1}{3}H_{3n}+H_{4n}\pmod{p}.\label{c-6}
\end{align}

Finally, noting $2n=\lfloor p/3\rfloor,3n=\lfloor p/2\rfloor, 4n=\lfloor 2p/3\rfloor$ and applying
\eqref{bL-1}--\eqref{bL-3} to the right-hand side of \eqref{c-6}, we arrive at the desired result:
\begin{align*}
\sum_{k=1}^{p-2}\sum_{j=1}^{p-2}\sum_{i=1}^{p-1}\frac{(-1)^k}{i(j+1)}{2j\choose j+k}{2i\choose i+k}
\equiv -q_p(3)\pmod{p},
\end{align*}
which is the case $p\equiv 1\pmod{3}$ of \eqref{c-1}.
\qed

\begin{lem}
For any prime $p\ge 5$, we have
\begin{align}
\sum_{k=1}^{p-1}\sum_{i=1}^{p-1}\sum_{j=0}^{p-1}\frac{(-1)^k}{i}{2i\choose i+k}{2j\choose j+k}\equiv -\left(\frac{p}{3}\right)q_p(3)\pmod{p}.
\label{c-7}
\end{align}
\end{lem}
{\noindent\it Proof.}
By \eqref{b-14} and \eqref{b-15}, we have
\begin{align}
\sum_{k=1}^{p-1}(-1)^k\sum_{i=1}^{p-1}\frac{1}{i}{2i\choose i+k}\sum_{j=0}^{p-1}{2j\choose j+k}
\equiv \sum_{k=1}^{p-1}\frac{(-1)^k(\alpha(k)-1)}{k}\left(\frac{p-k}{3}\right)\pmod{p}.
\label{c-8}
\end{align}
We shall only prove the case $p\equiv 1\pmod{3}$, and the proof of the case $p\equiv 2\pmod{3}$ runs similarly and the details are omitted.

Suppose that $p\equiv 1\pmod{3}$. Let $n=\lfloor p/6\rfloor$. Then $n\equiv -1/6\pmod{p}$. It is easy to check that
\begin{align*}
(-1)^k(\alpha(k)-1)\left(\frac{p-k}{3}\right)=
\begin{cases}
1\quad&\text{for $k\equiv 0\pmod{6}$},\\
0\quad&\text{for $k\equiv 1,4\pmod{6}$},\\
-1\quad&\text{for $k\equiv 2\pmod{6}$},\\
3\quad&\text{for $k\equiv 3\pmod{6}$},\\
-3\quad&\text{for $k\equiv 5\pmod{6}$}.
\end{cases}
\end{align*}
It follows that
\begin{align}
&\sum_{k=1}^{p-1}\sum_{i=1}^{p-1}\sum_{j=0}^{p-1}\frac{(-1)^k}{i}{2i\choose i+k}{2j\choose j+k}\notag\\[7pt]
&\equiv\sum_{k=1}^n\frac{1}{6k}-\sum_{k=1}^n\frac{1}{6k-4}+3\sum_{k=1}^n\frac{1}{6k-3}-3\sum_{k=1}^n\frac{1}{6k-1}\notag\\[7pt]
&\equiv\frac{1}{6}\left(\sum_{k=1}^n\frac{1}{k}-\sum_{k=1}^n\frac{1}{k+4n}+3\sum_{k=1}^n\frac{1}{k+3n}
-3\sum_{k=1}^n\frac{1}{k+n}\right)\notag\\[7pt]
&=\frac{1}{6}\left(4H_{n}-3H_{2n}-3H_{3n}+4H_{4n}-H_{5n}\right)\pmod{p}.\label{c-9}
\end{align}

Finally, applying \eqref{bL-1}--\eqref{bL-3} to the right-hand side of \eqref{c-9} gives
\begin{align*}
\sum_{k=1}^{p-1}\sum_{i=1}^{p-1}\sum_{j=0}^{p-1}\frac{(-1)^k}{i}{2i\choose i+k}{2j\choose j+k}\equiv -q_p(3)\pmod{p},
\end{align*}
which is the case $p\equiv 1\pmod{3}$ of \eqref{c-7}.
\qed

\begin{lem}
For any prime $p\ge 5$, we have
\begin{align}
\sum_{k=1}^{p-1}\sum_{i=1}^{p-2}\sum_{j=0}^{p-1}\frac{(-1)^k}{i+1}{2i\choose i+k}{2j\choose j+k}\equiv -\frac{1}{3}\left(2+\left(\frac{p}{3}\right)q_p(3)\right)\pmod{p}.\label{c-10}
\end{align}
\end{lem}
{\noindent\it Proof.}
By \eqref{b-14} and \eqref{b-16}, we have
\begin{align}
&\sum_{k=1}^{p-1}(-1)^k\sum_{i=1}^{p-2}\frac{1}{i+1}{2i\choose i+k}\sum_{j=0}^{p-1}{2j\choose j+k}\notag\\[7pt]
&\equiv \sum_{k=1}^{p-1}\left(-k\sum_{i=1}^{k-1}
\frac{\beta(i)+2}{i}+\sum_{i=1}^{k-1}\beta(i)+2k
+\frac{3}{2}\left(\frac{p}{3}\right)-\frac{3}{2}\right)\left(\frac{p-k}{3}\right)\notag\\[7pt]
&=-\sum_{k=1}^{p-1}k\left(\frac{p-k}{3}\right)\sum_{i=1}^{k-1}
\frac{\beta(i)+2}{i}+\sum_{k=1}^{p-1}\left(\frac{p-k}{3}\right)\sum_{i=1}^{k-1}\beta(i)\notag\\[7pt]
&+2\sum_{k=1}^{p-1}k\left(\frac{p-k}{3}\right)+\left(\frac{3}{2}\left(\frac{p}{3}\right)
-\frac{3}{2}\right)\sum_{k=1}^{p-1}\left(\frac{p-k}{3}\right)\pmod{p}.\label{cnew-1}
\end{align}

Similarly, we only prove the case $p\equiv 1\pmod{3}$. Suppose that $p\equiv 1\pmod{3}$.
It is easy to check that
\begin{align*}
k\left(\frac{p-k}{3}\right)=\begin{cases}
k\quad&\text{for $k\equiv 0\pmod{3}$},\\
0\quad&\text{for $k\equiv 1\pmod{3}$},\\
-k\quad&\text{for $k\equiv 2\pmod{3}$},
\end{cases}
\end{align*}
and
\begin{align*}
\left(\frac{p-k}{3}\right)\sum_{i=1}^{k-1}\beta(i)=
\begin{cases}
1\quad&\text{for $k\equiv 0\pmod{6}$},\\
0\quad&\text{for $k\equiv 1,4,5\pmod{6}$},\\
2\quad&\text{for $k\equiv 2\pmod{6}$},\\
-3\quad&\text{for $k\equiv 3\pmod{6}$}.
\end{cases}
\end{align*}
It follows that
\begin{align}
\Psi(m)&=\sum_{k=1}^{m}k\left(\frac{p-k}{3}\right)\notag\\[7pt]
&=\sum_{k=1}^{\lfloor m/3\rfloor}3k-\sum_{k=1}^{\lfloor (m+1)/3\rfloor}(3k-1)\notag\\[7pt]
&=\left\lfloor \frac{m+1}{3}\right\rfloor-(m+1)[3\mid m-2],\label{c-11}
\end{align}
and
\begin{align}
\sum_{k=1}^{p-1}\left(\frac{p-k}{3}\right)\sum_{i=1}^{k-1}\beta(i)=0.\label{c-12}
\end{align}
From \eqref{c-11}, we deduce that
\begin{align}
\Psi(p-1)\equiv -\frac{1}{3}\pmod{p}.\label{c-13}
\end{align}

Furthermore, we have
\begin{align}
&\sum_{k=1}^{p-1}k\left(\frac{p-k}{3}\right)\sum_{i=1}^{k-1}\frac{\beta(i)+2}{i}\notag\\[7pt]
&=\sum_{i=1}^{p-2}\frac{\beta(i)+2}{i}\sum_{k=i+1}^{p-1}k\left(\frac{p-k}{3}\right)\notag\\[7pt]
&=\sum_{i=1}^{p-2}\frac{\beta(i)+2}{i}\left(\Psi(p-1)-\Psi(i)\right)\notag\\[7pt]
&\equiv -\sum_{i=1}^{p-2}\frac{\beta(i)+2}{i}\left(\frac{1}{3}+\Psi(i)\right)\pmod{p},\label{c-14}
\end{align}
where we have used \eqref{c-13} in the last step. It is easy to check that
\begin{align}
\left(\beta(i)+2\right)\left(\frac{1}{3}+\Psi(i)\right)=
\begin{cases}
\frac{i+1}{3}\quad&\text{for $i\equiv 0\pmod{6}$},\\
0\quad&\text{for $i\equiv 1\pmod{6}$},\\
-\frac{2i+1}{3}\quad&\text{for $i\equiv 2\pmod{6}$},\\
i+1\quad&\text{for $i\equiv 3\pmod{6}$},\\
\frac{4i}{3}\quad&\text{for $i\equiv 4\pmod{6}$},\\
-(2i+1)\quad&\text{for $i\equiv 5\pmod{6}$}.
\end{cases}\label{c-15}
\end{align}
Combining \eqref{c-14} and \eqref{c-15} yields
\begin{align}
&\sum_{k=1}^{p-1}k\left(\frac{p-k}{3}\right)\sum_{i=1}^{k-1}\frac{\beta(i)+2}{i}\notag\\[7pt]
&\equiv -\frac{1}{3}\sum_{i=1}^{n-1}\frac{6i+1}{6i}+\frac{1}{3}\sum_{i=1}^{n}\frac{12i-7}{6i-4}
-\sum_{i=1}^{n}\frac{6i-2}{6i-3}-\frac{1}{3}\sum_{i=1}^{n}4
+\sum_{i=1}^{n}\frac{12i-1}{6i-1}\notag\\[7pt]
&=-\frac{1}{3}\sum_{i=1}^{n}\frac{1}{6i}+\frac{1}{3}\sum_{i=1}^{n}\frac{1}{6i-4}
-\sum_{i=1}^{n}\frac{1}{6i-3}+\sum_{i=1}^{n}\frac{1}{6i-1}+\frac{1}{3}+\frac{1}{18n}\notag\\[7pt]
&\equiv -\frac{1}{18}\sum_{i=1}^{n}\frac{1}{i}+\frac{1}{18}\sum_{i=1}^{n}\frac{1}{i+4n}
-\frac{1}{6}\sum_{i=1}^{n}\frac{1}{i+3n}+\frac{1}{6}\sum_{i=1}^{n}\frac{1}{i+n}\notag\\[7pt]
&=\frac{1}{6}\left(-\frac{4}{3}H_n+H_{2n}+H_{3n}-\frac{4}{3}H_{4n}+\frac{1}{3}H_{5n}\right)\notag\\[7pt]
&\equiv \frac{1}{3}q_p(3)\pmod{p},\label{c-16}
\end{align}
where we have used \eqref{bL-1}--\eqref{bL-3} in the last step.

Finally, substituting \eqref{c-12}, \eqref{c-13} and \eqref{c-16} into the right-hand side of \eqref{cnew-1} gives
\begin{align*}
\sum_{k=1}^{p-1}\sum_{i=1}^{p-2}\sum_{j=0}^{p-1}\frac{(-1)^k}{i+1}{2i\choose i+k}{2j\choose j+k}
\equiv -\frac{1}{3}\left(q_p(3)+2\right)\pmod{p},
\end{align*}
which is the case $p\equiv 1\pmod{3}$ of \eqref{c-10}.
\qed

\section{Proof of \eqref{a-2} and \eqref{a-3}}
The proof of \eqref{a-3} follows from \eqref{a-1} and \eqref{a-2}, and it remains to prove \eqref{a-2}.

By \eqref{a-5}, we have
\begin{align}
\sum_{k=0}^{p-1}M_k^2&=\sum_{k=0}^{p-1}\sum_{i=0}^k \sum_{j=0}^k(-1)^{i+j}{k\choose i}{k\choose j}C_{i+1}C_{j+1}\notag\\
&=\sum_{i=0}^{p-1} \sum_{j=0}^{p-1}(-1)^{i+j}C_{i+1}C_{j+1}\sum_{k=0}^{p-1}{k\choose i}{k\choose j}.
\label{d-1}
\end{align}

From the identity \cite[(9), page 15]{riordan-b-1968}:
\begin{align*}
{k\choose i}{k\choose j}=\sum_{l=0}^j
{l+i\choose j}{j\choose l}{k\choose l+i},
\end{align*}
we deduce that
\begin{align}
\sum_{k=0}^{p-1}{k\choose i}{k\choose j}
&=\sum_{l=0}^j
{l+i\choose j}{j\choose l}\sum_{k=0}^{p-1}{k\choose l+i}\notag\\
&=\sum_{l=0}^j
{i+l\choose j}{j\choose l}{p\choose i+l+1},\label{d-2}
\end{align}
where we have utilized the induction on $p$ in the last step.
Letting $l\to k-i$ on the right-hand side of \eqref{d-2}, we rewrite \eqref{d-2} as
\begin{align*}
\sum_{k=0}^{p-1}{k\choose i}{k\choose j}=p\sum_{k=0}^{i+j}
\frac{1}{k+1}{k\choose j}{j\choose k-i}{p-1\choose k}.
\end{align*}

Note that ${k\choose j}{j\choose k-i}\equiv 0\pmod{p}$ for $0\le i\le p-1,0\le j\le p-1$ and $p\le k \le 2p-2$, and
\begin{align*}
\frac{p}{k+1}{p-1\choose k}\equiv \frac{p(-1)^k}{k+1}\pmod{p^2},
\end{align*}
for $0\le k\le p-1$. It follows that
\begin{align}
\sum_{k=0}^{p-1}{k\choose i}{k\choose j}\equiv p\sum_{k=0}^{i+j}
\frac{(-1)^k}{k+1}{k\choose j}{j\choose k-i}\pmod{p^2}.\label{d-3}
\end{align}

Recall the following partial fraction decomposition:
\begin{align}
\sum_{k=0}^{i+j}\frac{(-1)^{i+j+k}}{x+k}{k\choose j}{j\choose k-i}=
\frac{(x)_i (x)_j}{(x)_{i+j+1}},\label{d-4}
\end{align}
where $(x)_0=1$ and $(x)_k=x(x+1)\cdots (x+k-1)$ for $k\ge 1$.
Letting $x=1$ in \eqref{d-4} gives
\begin{align}
\sum_{k=0}^{i+j}\frac{(-1)^{k}}{k+1}{k\choose j}{j\choose k-i}=\frac{(-1)^{i+j}}{(i+j+1){i+j\choose j}}.\label{d-5}
\end{align}
It follows from \eqref{d-3} and \eqref{d-5} that
\begin{align}
\sum_{k=0}^{p-1}{k\choose i}{k\choose j}\equiv  \frac{p(-1)^{i+j}}{(i+j+1){i+j\choose j}}\pmod{p^2}.\label{dnew-1}
\end{align}
Combining \eqref{d-1} and \eqref{dnew-1} gives
\begin{align}
\sum_{k=0}^{p-1}M_k^2
\equiv p\sum_{i=0}^{p-1} \sum_{j=0}^{p-1}\frac{C_{i+1}C_{j+1}}{(i+j+1){i+j\choose j}}\pmod{p^2}.
\label{d-6}
\end{align}

Let
\begin{align*}
S(i,j)=\frac{{2i\choose i}{2j\choose j}}{{i+j\choose i}}.
\end{align*}
From the identity due to Von Szily \cite{szily-1984} (see also \cite[(3.38)]{gould-b-1972}):
\begin{align}
S(i,j)=\sum_{k}(-1)^k{2i\choose i+k}{2j\choose j+k},\label{d-7}
\end{align}
we deduce that $S(i,j)$ is an integer. It is easy to verify that the numbers $S(i,j)$ satisfy the recurrence:
\begin{align}
4S(i,j)=S(i+1,j)+S(i,j+1).\label{d-8}
\end{align}

Note that
\begin{align}
\frac{C_{i+1}C_{j+1}}{(i+j+1){i+j\choose j}}=\frac{(i+j+2)S(i+1,j+1)}{(i+1)(j+1)(i+2)(j+2)}.
\label{d-9}
\end{align}
It follows from \eqref{d-6} and \eqref{d-9} that
\begin{align}
\sum_{k=0}^{p-1}M_k^2
&\equiv p\sum_{i=0}^{p-1} \sum_{j=0}^{p-1}\frac{(i+j+2)S(i+1,j+1)}{(i+1)(j+1)(i+2)(j+2)}\notag\\
&=p\sum_{i=1}^{p} \sum_{j=1}^{p}\frac{(i+j)S(i,j)}{ij(i+1)(j+1)}\notag\\
&=2p\sum_{i=1}^{p} \sum_{j=1}^{p}\frac{S(i,j)}{i(i+1)(j+1)}\notag\\
&=2p\sum_{i=1}^{p} \sum_{j=1}^{p}\frac{S(i,j)}{i(j+1)}-2p\sum_{i=1}^{p} \sum_{j=1}^{p}\frac{S(i,j)}{(i+1)(j+1)}\pmod{p^2},\label{d-10}
\end{align}
where we have used the symmetry with respect to $i$ and $j$ in the third step.

By \eqref{d-8}, we have
\begin{align}
&2p\sum_{i=1}^{p}\sum_{j=1}^p\frac{S(i,j)}{(i+1)(j+1)}\notag\\
&=\frac{p}{2}\sum_{i=1}^{p}\sum_{j=1}^p\frac{S(i+1,j)+S(i,j+1)}{(i+1)(j+1)}\notag\\
&=p\sum_{i=1}^{p}\sum_{j=1}^p\frac{S(i+1,j)}{(i+1)(j+1)}\notag\\
&=p\sum_{i=2}^{p+1}\sum_{j=1}^p\frac{S(i,j)}{i(j+1)}\notag\\
&=p\sum_{i=1}^{p}\sum_{j=1}^p\frac{S(i,j)}{i(j+1)}+p\sum_{j=1}^p\frac{S(p+1,j)}{(p+1)(j+1)}
-p\sum_{j=1}^p\frac{S(1,j)}{j+1}.\label{d-11}
\end{align}

Furthermore, by \eqref{b-2} we have
\begin{align*}
p\sum_{j=1}^p\frac{S(p+1,j)}{(p+1)(j+1)}&=\frac{2p(2p+1)}{(p+1)^2}{2p\choose p}\sum_{j=1}^p\frac{{2j\choose j}}{(j+1){p+1+j\choose j}}\\
&\equiv \frac{4p(2p+1)}{(p+1)^2}\sum_{j=1}^p\frac{{2j\choose j}}{(j+1){p+1+j\choose j}}\pmod{p^2}.
\end{align*}
For $1\le j\le p-2$, we have ${p+1+j\choose j}\equiv j+1\pmod{p}$. It follows that
\begin{align}
p\sum_{j=1}^p\frac{S(p+1,j)}{(p+1)(j+1)}&\equiv
4p\sum_{j=1}^{p-2}\frac{{2j\choose j}}{(j+1)^2}+\frac{2(2p+1)}{(p+1)(2p-1)}
+\frac{4p}{(p+1)^2}\notag\\
&\equiv 12p\left(\frac{p}{3}\right)-2-2p\pmod{p^2},\label{d-12}
\end{align}
where we have used \eqref{b-18} in the last step.

On the other hand, we have
\begin{align}
p\sum_{j=1}^p\frac{S(1,j)}{j+1}&=2p\sum_{j=1}^p\frac{{2j\choose j}}{(j+1)^2}\notag\\
&=2p\sum_{j=1}^{p-2}\frac{{2j\choose j}}{(j+1)^2}+\frac{2}{2p-1}{2p-1\choose p-1}+\frac{2p}{(p+1)^2}{2p\choose p}\notag\\
&\equiv 6p\left(\frac{p}{3}\right)-2\pmod{p^2},\label{d-13}
\end{align}
where we have used \eqref{b-1}, \eqref{b-2} and \eqref{b-18} in the last step.

It follows from \eqref{d-11}, \eqref{d-12} and \eqref{d-13} that
\begin{align}
2p\sum_{i=1}^{p}\sum_{j=1}^p\frac{S(i,j)}{(i+1)(j+1)}-p\sum_{i=1}^{p}\sum_{j=1}^p\frac{S(i,j)}{i(j+1)}
\equiv 2p\left(3\left(\frac{p}{3}\right)-1\right)\pmod{p^2}.\label{d-14}
\end{align}
Combining \eqref{d-10} and \eqref{d-14} gives
\begin{align}
\sum_{k=0}^{p-1}M_k^2\equiv p\sum_{i=1}^{p}\sum_{j=1}^p\frac{S(i,j)}{i(j+1)}
-2p\left(3\left(\frac{p}{3}\right)-1\right)\pmod{p^2}.\label{d-15}
\end{align}

Note that
\begin{align}
&p\sum_{i=1}^{p}\sum_{j=1}^p\frac{S(i,j)}{i(j+1)}\notag\\
&=p\sum_{i=1}^{p-1}\sum_{j=1}^p\frac{S(i,j)}{i(j+1)}+\sum_{j=1}^p\frac{S(p,j)}{j+1}\notag\\
&=p\sum_{i=1}^{p-1}\sum_{j=1}^{p-2}\frac{S(i,j)}{i(j+1)}+\sum_{j=1}^p\frac{S(p,j)}{j+1}
+\sum_{i=1}^{p-1}\frac{S(i,p-1)}{i}+p\sum_{i=1}^{p-1}\frac{S(i,p)}{i(p+1)}.\label{d-16}
\end{align}

Furthermore, we have
\begin{align}
p\sum_{i=1}^{p-1}\frac{S(i,p)}{i(p+1)}\equiv 2p\sum_{i=1}^{p-1}\frac{{2i\choose i}}{i}\equiv 0\pmod{p^2},
\label{d-17}
\end{align}
where we have used \eqref{b-2}, \eqref{bnew-4} and the fact that ${p+i\choose i}\equiv 1\pmod{p}$ for
$1\le i\le p-1$.

For $1\le i\le p-1$, we have $p/i{i+p-1\choose i}\equiv 1-pH_{i-1}\pmod{p^2}$, and so
\begin{align}
\sum_{i=1}^{p-1}\frac{S(i,p-1)}{i}
&=\frac{p}{2p-1}{2p-1\choose p-1}\sum_{i=1}^{p-1}\frac{{2i\choose i}}{i{i+p-1\choose i}}\notag\\[7pt]
&\equiv -(2p+1)\sum_{i=1}^{p-1}{2i\choose i}+p\sum_{i=1}^{p-1}{2i\choose i}H_{i-1}\notag\\[7pt]
&=-(2p+1)\sum_{i=1}^{p-1}{2i\choose i}+p\sum_{i=1}^{p-1}{2i\choose i}H_{i}
-p\sum_{i=1}^{p-1}\frac{{2i\choose i}}{i}\notag\\[7pt]
&\equiv 2p+1-\left(\frac{p}{3}\right)\left(2p+3^{p-1}\right)\pmod{p^2},\label{d-18}
\end{align}
where we have used \eqref{b-1}, \eqref{b-3}, \eqref{bnew-4} and \eqref{b-8}.

For $1\le j\le p-1$, we have ${p+j\choose j}\equiv 1+pH_j \pmod{p^2}$, and so
\begin{align}
\sum_{j=1}^p\frac{S(p,j)}{j+1}&=\frac{{2p\choose p}}{p+1}
+{2p\choose p}\sum_{j=1}^{p-1}\frac{{2j\choose j}}{(j+1){p+j\choose j}}\notag\\[7pt]
&\equiv 2-2p+2\sum_{j=1}^{p-1}C_j(1-pH_j)\notag\\[7pt]
&\equiv 3^p\left(\frac{p}{3}\right)-2p-1\pmod{p^2},\label{d-19}
\end{align}
where we have used \eqref{b-2}, \eqref{b-4} and \eqref{b-9}.

It follows from \eqref{d-16}--\eqref{d-19} that
\begin{align}
p\sum_{i=1}^{p}\sum_{j=1}^p\frac{S(i,j)}{i(j+1)}
\equiv p\sum_{i=1}^{p-1}\sum_{j=1}^{p-2}\frac{S(i,j)}{i(j+1)}
+2\left(\frac{p}{3}\right)\left(3^{p-1}-p\right)\pmod{p^2}.\label{d-20}
\end{align}
Combining \eqref{d-15} and \eqref{d-20} gives
\begin{align}
\sum_{k=0}^{p-1}M_k^2\equiv p\sum_{i=1}^{p-1}\sum_{j=1}^{p-2}\frac{S(i,j)}{i(j+1)}
+2\left(\frac{p}{3}\right)\left(3^{p-1}-4p\right)+2p\pmod{p^2}.\label{d-21}
\end{align}

By \eqref{d-7}, we have
\begin{align}
S(i,j)=2\sum_{k=1}^{j}(-1)^k{2i\choose i+k}{2j\choose j+k}+{2i\choose i}{2j\choose j}.\label{dnew-2}
\end{align}
It follows that
\begin{align}
&p\sum_{i=1}^{p-1}\sum_{j=1}^{p-2}\frac{S(i,j)}{i(j+1)}\notag\\[7pt]
&=2p\sum_{i=1}^{p-1}\sum_{j=1}^{p-2}\sum_{k=1}^{j}\frac{(-1)^k}{i(j+1)}{2i\choose i+k}{2j\choose j+k}
+p\sum_{i=1}^{p-1}\sum_{j=1}^{p-2}\frac{{2i\choose i}{2j\choose j}}{i(j+1)}\notag\\[7pt]
&\equiv 2p\sum_{k=1}^{p-2}\sum_{i=1}^{p-1}\sum_{j=1}^{p-2}\frac{(-1)^k}{i(j+1)}{2i\choose i+k}{2j\choose j+k}\pmod{p^2},\label{d-22}
\end{align}
where we have used \eqref{bnew-4} in the last step.

Finally, combining \eqref{d-21} and \eqref{d-22} yields
\begin{align}
\sum_{k=0}^{p-1}M_k^2&\equiv 2p\sum_{k=1}^{p-2}\sum_{i=1}^{p-1}\sum_{j=1}^{p-2}\frac{(-1)^k}{i(j+1)}{2i\choose i+k}{2j\choose j+k}\notag\\[7pt]
&+2\left(\frac{p}{3}\right)\left(3^{p-1}-4p\right)+2p\pmod{p^2}.\label{d-23}
\end{align}
Then the proof of \eqref{a-2} follows from \eqref{c-1} and \eqref{d-23}.

\section{Proof of \eqref{a-4}}
By \eqref{a-5} and \eqref{a-6}, we have
\begin{align}
\sum_{k=0}^{p-1}T_kM_k&=\sum_{k=0}^{p-1}
\sum_{i=0}^k \sum_{j=0}^k(-1)^{i+j}{k\choose i}{k\choose j}{2j\choose j}C_{i+1}\notag\\[7pt]
&=\sum_{i=0}^{p-1} \sum_{j=0}^{p-1}(-1)^{i+j}{2j\choose j}C_{i+1}\sum_{k=0}^{p-1}{k\choose i}{k\choose j}.\label{e-1}
\end{align}
Applying \eqref{dnew-1} to the right-hand side of \eqref{e-1}, we obtain
\begin{align*}
\sum_{k=0}^{p-1}T_kM_k\equiv p\sum_{i=0}^{p-1} \sum_{j=0}^{p-1}\frac{{2j\choose j}C_{i+1}}{(i+j+1){i+j\choose j}}\pmod{p^2}.
\end{align*}
Noting that
\begin{align*}
\frac{{2j\choose j}C_{i+1}}{(i+j+1){i+j\choose j}}=\frac{S(i+1,j)}{(i+1)(i+2)},
\end{align*}
we have
\begin{align}
\sum_{k=0}^{p-1}T_kM_k&\equiv p\sum_{i=0}^{p-1} \sum_{j=0}^{p-1}\frac{S(i+1,j)}{(i+1)(i+2)}\notag\\[7pt]
&=p\sum_{i=1}^{p} \sum_{j=0}^{p-1}\frac{S(i,j)}{i(i+1)}\notag\\[7pt]
&\equiv p\sum_{i=1}^p\sum_{j=0}^{p-1}\frac{S(i,j)}{i}-p\sum_{i=1}^p\sum_{j=0}^{p-1}\frac{S(i,j)}{i+1}\pmod{p^2}.
\label{e-2}
\end{align}

Note that
\begin{align}
p\sum_{i=1}^p\sum_{j=0}^{p-1}\frac{S(i,j)}{i}=p\sum_{i=1}^{p-1}\sum_{j=0}^{p-1}\frac{S(i,j)}{i}
+\sum_{j=0}^{p-1}S(p,j),\label{e-3}
\end{align}
and
\begin{align}
p\sum_{i=1}^p\sum_{j=0}^{p-1}\frac{S(i,j)}{i+1}=p\sum_{i=1}^{p-2}\sum_{j=0}^{p-1}\frac{S(i,j)}{i+1}
+\sum_{j=0}^{p-1}S(p-1,j)+p\sum_{j=0}^{p-1}\frac{S(p,j)}{p+1}.\label{e-4}
\end{align}

For $0\le j\le p-1$, we have ${p+j\choose j}\equiv 1+pH_j\pmod{p^2}$, and so
\begin{align}
\sum_{j=0}^{p-1}S(p,j)&=
{2p\choose p}\sum_{j=0}^{p-1}\frac{{2j\choose j}}{{p+j\choose j}}\notag\\
&\equiv 2\sum_{j=0}^{p-1}{2j\choose j}(1-pH_j)\notag\\
&\equiv 2\left(\frac{p}{3}\right)3^{p-1}\pmod{p^2},\label{e-5}
\end{align}
where we have used \eqref{b-2}, \eqref{b-3} and \eqref{b-8}.

By \eqref{b-2} and \eqref{b-3}, we have
\begin{align}
p\sum_{j=0}^{p-1}\frac{S(p,j)}{p+1}\equiv p\sum_{j=0}^{p-1}S(p,j)\equiv 2p\left(\frac{p}{3}\right)\pmod{p^2}.\label{e-6}
\end{align}

By \eqref{b-1} and $p/{p+j-1\choose j}\equiv j(1-pH_{j-1})\pmod{p^2}$ for $0\le j\le p-1$, we have
\begin{align}
\sum_{j=0}^{p-1}S(p-1,j)&=\frac{p}{2p-1}{2p-1\choose p-1}\sum_{j=0}^{p-1}\frac{{2j\choose j}}{{j+p-1\choose j}}\notag\\[7pt]
&\equiv -(2p+1)\sum_{j=0}^{p-1}j{2j\choose j}+p\sum_{j=0}^{p-1}j{2j\choose j}H_{j-1}\notag\\[7pt]
&=-(2p+1)\sum_{j=0}^{p-1}j{2j\choose j}+p\sum_{j=0}^{p-1}j{2j\choose j}H_{j}-p\sum_{j=0}^{p-1}{2j\choose j}\notag\\[7pt]
&\equiv -\frac{p}{3}+\frac{2}{3}\left(\frac{p}{3}\right)3^{p-1}\pmod{p^2},\label{e-7}
\end{align}
where we have used \eqref{b-3}, \eqref{b-5} and \eqref{bnew-3} in the last step.

It follows from \eqref{e-2}--\eqref{e-7} that
\begin{align}
\sum_{k=0}^{p-1}T_kM_k&\equiv
p\sum_{i=1}^{p-1}\sum_{j=0}^{p-1}\frac{S(i,j)}{i}-p\sum_{i=1}^{p-2}\sum_{j=0}^{p-1}\frac{S(i,j)}{i+1}\notag\\[7pt]
&+\frac{p}{3}-2p\left(\frac{p}{3}\right)+\frac{4}{3}\left(\frac{p}{3}\right)3^{p-1}\pmod{p^2}.\label{e-8}
\end{align}

By \eqref{dnew-2}, we have
\begin{align}
&p\sum_{i=1}^{p-1}\sum_{j=0}^{p-1}\frac{S(i,j)}{i}\notag\\[7pt]
&=2p\sum_{i=1}^{p-1}\sum_{j=0}^{p-1}\sum_{k=1}^{j}\frac{(-1)^k}{i}{2i\choose i+k}{2j\choose j+k} +p\sum_{i=1}^{p-1}\sum_{j=0}^{p-1}\frac{{2i\choose i}{2j\choose j}}{i}\notag\\[7pt]
&\equiv 2p\sum_{k=1}^{p-1}\sum_{i=1}^{p-1}\sum_{j=0}^{p-1}\frac{(-1)^k}{i}{2i\choose i+k}{2j\choose j+k}\pmod{p^2},\label{e-9}
\end{align}
where we have used \eqref{bnew-4} in the last step.

On the other hand, we have
\begin{align}
&p\sum_{i=1}^{p-2}\sum_{j=0}^{p-1}\frac{S(i,j)}{i+1}\notag\\[7pt]
&=2p\sum_{i=1}^{p-2}\sum_{j=0}^{p-1}\sum_{k=1}^j\frac{(-1)^k}{i+1}{2i\choose i+k}{2j\choose j+k}
 +p\sum_{i=1}^{p-2}\sum_{j=0}^{p-1}\frac{{2i\choose i}{2j\choose j}}{i+1}\notag\\[7pt]
&\equiv 2p\sum_{k=1}^{p-1}\sum_{i=1}^{p-2}\sum_{j=0}^{p-1}\frac{(-1)^k}{i+1}{2i\choose i+k}{2j\choose j+k}+\frac{p}{2}\left(3- \left(\frac{p}{3}\right)\right)\pmod{p^2},\label{e-10}
\end{align}
where we have used \eqref{b-3}, \eqref{b-4} and $C_{p-1}\equiv -1\pmod{p}$ in the last step.

It follows from \eqref{e-8}--\eqref{e-10} that
\begin{align}
&\sum_{k=0}^{p-1}T_kM_k\notag\\[7pt]
&\equiv2p\sum_{k=1}^{p-1}\sum_{i=1}^{p-1}\sum_{j=0}^{p-1}\frac{(-1)^k}{i}{2i\choose i+k}{2j\choose j+k}\notag\\[7pt]
&-2p\sum_{k=1}^{p-1}\sum_{i=1}^{p-2}\sum_{j=0}^{p-1}\frac{(-1)^k}{i+1}{2i\choose i+k}{2j\choose j+k}\notag\\[7pt]
&-\frac{7p}{6}-\frac{3p}{2}\left(\frac{p}{3}\right)+\frac{4}{3}\left(\frac{p}{3}\right)3^{p-1}\pmod{p^2}.
\label{e-11}
\end{align}
Then the proof of \eqref{a-4} follows from \eqref{c-7}, \eqref{c-10} and \eqref{e-11}.

\vskip 5mm \noindent{\bf Acknowledgments.}
This work was supported by the National Natural Science Foundation of China (grant 12171370).

\end{document}